\documentclass[12pt]{article}
\author{David J. Platt\footnote{Supported by ARC Discovery Project DP160100932 and  EPSRC Grant EP/K034383/1.}\\ School of Mathematics, \\ University of Bristol, Bristol, UK\\dave.platt@bris.ac.uk\and  Timothy S. Trudgian\footnote{Supported by ARC Discovery Project DP160100932 and ARC Future Fellowship FT160100094.} \\
School of Science\\ The University of New South Wales Canberra, Australia \\
  t.trudgian@adfa.edu.au}
  \title{The error term in the prime number theorem}
\usepackage{indentfirst}
\usepackage{url}
\usepackage{enumerate}
\usepackage{amsthm}
\usepackage{amsmath}
\usepackage{comment}
\usepackage{fullpage}
\usepackage{amssymb}
\usepackage{booktabs}

\newtheorem{thm}{Theorem}

\newtheorem{Lem}{Lemma}

\newtheorem{cor}{Corollary}

\newcommand{\dif}{\textrm d}
\begin{document}
\maketitle
\begin{abstract}
\noindent 
We make explicit a theorem of Pintz, which gives a version of the prime number theorem with error term roughly square-root of that which was previously known. We apply this to a long-standing problem concerning an inequality studied by Ramanujan.
\end{abstract}
\section{Introduction}
\noindent
The zero-free region of the Riemann zeta-function $\zeta(s)$ is intimately connected with the size of the error term in the prime number theorem $\psi(x) \sim x$ where $\psi(x) = \sum_{p^{m} \leq x} \log p$. Ingham \cite{Ingham} essentially showed that if $\zeta(s)$ has no zeroes with real part $\sigma \geq 1 - \eta(t)$, where $\eta(t)$ is a decreasing function, then one has
\begin{equation}\label{spring}
\frac{\psi(x) -x}{x} \ll  \exp\left\{ -\frac{1}{2}  (1-\epsilon) \omega(x)\right\},
\end{equation}
where
\begin{equation}\label{rugby}
\omega(x):= \min_{t\geq 1} \{\eta(t) \log x + \log t\}.
\end{equation}
The classical zero-free region has $\eta(t) = 1/(R\log t)$, for some positive constant $R$. This gives $t = \exp\{\sqrt{\log x /R}\}$ as the optimal choice for $t$ in (\ref{rugby}) and gives
\begin{equation}\label{tour}
\frac{\psi(x) -x}{x} \ll \exp\left( -(1-\epsilon)\sqrt{\frac{\log x}{R}}\right).
\end{equation}
Explicit versions of (\ref{tour}) have been given in \cite{Dusart, Rosser, Schoenfeld, Trudgian}. These have the form
\begin{equation*}\label{TV}
\left|\frac{\psi(x) -x}{x}\right| \leq A (\log x)^{B} \exp\left( -\sqrt{\frac{\log x}{R}}\right), \quad (x\geq x_{0}),
\end{equation*}
with specific values of $A, B$ and $x_{0}$. The $(\log x)^{B}$ factor is not very important, but does give an explicit version of the $\epsilon$ in (\ref{tour}). 
%All explicit estimates have focussed on the constant $A$, and slight improvements on $R$.

Pintz \cite{Pintz} showed that Ingham's bound is a substantial overestimate: we can delete the factor of $1/2$ in (\ref{spring}), which leads to replacing the $1-\epsilon$ in (\ref{tour}) by $2-\epsilon$. Therefore, making Pintz's theorem explicit obtains an error term in the prime number theorem that is almost the square-root of the current bound. We prove this in the following theorem.

\begin{thm}\label{caprice}
Let $R = 5.573412$. For each row $\{X, A, B, C, \epsilon_0\}$ from Table \ref{tab3} we have
\begin{equation}\label{marcellina}
\left|\frac{\psi(x) - x}{x}\right| \leq A \left(\frac{\log x}{R}\right)^B \exp\left(-C\sqrt{\frac{\log x}{R}}\right)
\end{equation}
and
\begin{equation*}
\left|\psi(x)-x\right|\leq \epsilon_0 x
\end{equation*}
for all $\log x \geq X$.
\begin{comment}
Further, as $x\rightarrow\infty$ we have $B,C\rightarrow 2$ and $A\rightarrow 192$. Moreover, (\ref{marcellina}) holds for any $R\leq 5.573412$ such that $\zeta(s)$ has no zeroes for $\sigma \geq 1 - (R\log t)^{-1}$ for $t\geq 3$.
\end{comment}
\end{thm}

We also state the following obvious corollary for $\theta(x) = \sum_{p \leq x} \log p$.
\begin{cor}\label{caprice1}
For each row $\{X, A, B, C\}$ from Table \ref{tab3} we have
$$
\left|\frac{\theta(x) - x}{x}\right| \leq A_1 \left(\frac{\log x}{R}\right)^B \exp\left(-C\sqrt{\frac{\log x}{R}}\right)\;\;\;\textrm{for all } \log x\geq X,
$$
where $A_1 = A + 0.1$.
\begin{proof}
This follows trivially (and wastefully) from the work of Dusart  \cite[Cor.\ 4.5]{Dusart} or the authors \cite[Cor.\ 2]{PT}.
\begin{equation}\label{shoe}
\psi(x)-\theta(x)<(1+1.47\cdot 10^{-7})\sqrt{x}+1.78x^{1/3}.
\end{equation}
\end{proof}
\end{cor}
We remark that the values of $\epsilon_{0}$ in Table \ref{tab3} improve on those given by Faber and Kadiri in \cite{Faber_fix} when $\log x \geq 3000$. See also Broadbent et al.\ \cite{Broad} for some recent results, and an excellent survey of the field.

We also give an explicit estimate of another version of the prime number theorem.
\begin{cor}\label{white}
%\begin{comment}
For all $\log x\geq 2000$ we have
\begin{equation}\label{library}
|\pi(x) - \emph{\textrm{li}}(x)| \leq 235x(\log x)^{0.52} \exp(-0.8\sqrt{\log x}).
\end{equation}
%\end{comment}
\end{cor}

This improves on results by Dusart \cite[Thm 1.12]{DusartT} and by Trudgian \cite[Thm 2]{Trudgian} for all $\log x \geq 2000$. We note that, with a little work, the coefficient $0.8$ in (\ref{library}) could be increased but at present not beyond $0.8471\ldots$.

\begin{table}[!h]
\centering
\caption{Values of $X, A, B, C, \epsilon_0$ in Theorem \ref{caprice}.}
\label{tab3}
\begin{tabular}{cccccc}
$X$ & $\sigma$ & $A$ & $B$ & $C$ & $\epsilon_0$\\
$1\,000$ & $0.98$ &  $461.9$ & $1.52$ & $1.89$ & $1.20\cdot 10^{-5}$\\
$2\,000$ & $0.98$ &  $411.4$ & $1.52$ & $1.89$ & $8.35\cdot 10^{-10}$\\
$3\,000$ & $0.98$ &  $379.6$ & $1.52$ & $1.89$ & $4.51\cdot 10^{-13}$\\
$4\,000$ & $0.98$ &  $356.3$ & $1.52$ & $1.89$ & $7.33\cdot 10^{-16}$\\
$5\,000$ & $0.99$ &  $713.0$ & $1.51$ & $1.94$ & $9.77\cdot 10^{-19}$\\
$6\,000$ & $0.99$ &  $611.6$ & $1.51$ & $1.94$ & $4.23\cdot 10^{-21}$\\
$7\,000$ & $0.99$ &  $590.1$ & $1.51$ & $1.94$ & $3.09\cdot 10^{-23}$\\
$8\,000$ & $0.99$ &  $570.5$ & $1.51$ & $1.94$ & $3.12\cdot 10^{-25}$\\
$9\,000$ & $0.99$ &  $552.3$ & $1.51$ & $1.94$ & $4.11\cdot 10^{-27}$\\
$10\,000$ & $0.99$ &  $535.4$ & $1.51$ & $1.94$ & $6.78\cdot 10^{-29}$
%$1\,000$ & $0.98$ & $273.2$ & $1.52$ & $1.89$ & $7.06\cdot 10^{-6}$\\
%$2\,000$ & $0.98$ & $273.0$ & $1.52$ & $1.89$ & $5.54\cdot 10^{-10}$\\
%$2\,750$ & $0.99$ & $944.9$ & $1.51$ & $1.94$ & $1.84\cdot 10^{-12}$\\
%$3\,000$ & $0.99$ & $637.7$ & $1.51$ & $1.94$ & $2.07\cdot 10^{-13}$\\
%$4\,000$ & $0.99$ & $296.7$ & $1.51$ & $1.94$ & $1.37\cdot 10^{-16}$\\
%$5\,000$ & $0.99$ & $281.8$ & $1.51$ & $1.94$ & $3.87\cdot 10^{-19}$\\
%$6\,000$ & $0.99$ & $281.5$ & $1.51$ & $1.94$ & $1.95\cdot 10^{-21}$\\
%$7\,000$ & $0.99$ & $281.4$ & $1.51$ & $1.94$ & $1.47\cdot 10^{-23}$\\
%$8\,000$ & $0.99$ & $281.3$ & $1.51$ & $1.94$ & $1.54\cdot 10^{-25}$\\
%$9\,000$ & $0.99$ & $281.3$ & $1.51$ & $1.94$ & $2.09\cdot 10^{-27}$\\
%$10\,000$ & $0.99$ & $281.2$ & $1.51$ & $1.94$ & $3.56\cdot 10^{-29}$
\end{tabular}
\end{table}

We collect some lemmas in \S \ref{prep} that allow us to prove Theorem \ref{caprice} in \S \ref{whim}. We then prove Corollary \ref{white} in \S \ref{hunger}. While there are many applications for Theorem \ref{caprice} we focus on just one in \S \ref{bartolo}, where we prove the following theorem.
\begin{thm}\label{sidelines}
%\begin{comment}
The inequality
\begin{equation*}\label{oranges}
\pi^{2} (x) < \frac{e x}{\log x} \pi \left(\frac{x}{e}\right)
\end{equation*}
is true for all $38\,358\,837\,682 < x \leq \exp(58)$ and for $x\geq \exp(3\,915)$.
%\end{comment}
\end{thm}
Throughout this paper, we use the notation $f(x) = \vartheta(g(x))$ to mean that $|f(x)| \leq g(x)$.

%We note that we can convert the bounds on $\psi(x) -x$ into bounds for $\theta(x) -x$ or bounds for $\pi(x) - \textrm{li}(x)$ without too much trouble.

\section{Preparatory lemmas}\label{prep}

We start with an explicit version of the explicit formula, given in \cite[Thm.\ 1.3]{Dudek}.

\begin{Lem}
Let $50<T<x$ where $x> e^{60}$ and $x$ is half an odd integer. Then
\begin{equation}\label{honduras} 
\frac{\psi(x)-x}{x} = -\sum_{|\gamma| <T} \frac{x^{\rho -1}}{\rho} + \vartheta\left( \frac{2 \log^{2} x}{T}\right).
\end{equation}
\end{Lem}

We also need a result on sums over zeroes of $\zeta(s)$, which we quote from \cite[Lem.\ 2.10]{Saouter}.
%[Demichel, Saouter, and Trudgian]
\begin{Lem} \label{happy}
If $T\geq 2\pi e$, then
\begin{equation*}\label{snooker}
\sum_{0< \gamma \leq T} \frac{1}{\gamma} = \frac{1}{4\pi} \left( \log \frac{T}{2\pi}\right)^{2} + \vartheta(0.9321).
\end{equation*}
\end{Lem}

We shall make use of the following zero-free result  \cite{Hoff}.
\begin{Lem}\label{sneezy}
There are no zeroes of $\zeta(s)$ in the region $\sigma \geq 1 - 1/(R \log t)$ for $t\geq 3$ where $R = 5.573412$.
\end{Lem}
We remark that using Ford's   result \cite{Ford} that $\zeta(s)$ has no zeroes in the region $\sigma \geq 1 - 1/(57.54 (\log t)^{2/3} (\log\log t)^{1/3})$ for $t\geq 3$,  would give a totally explicit version of the estimate
$\psi(x)- x = O(x\exp(-c (\log x)^{3/5} (\log\log x)^{-1/5}))$, which is the best asymptotic estimate for the error term in the prime number theorem. We do not pursue this here.

We also need the height to which the Riemann hypothesis has been verified.

\begin{Lem}\label{pen}
Let $0< \beta <1$. Then, if $\zeta(\beta + i\gamma) =0$ and $|\gamma| \leq 3\cdot 10^{12}$ we have $\beta = \frac{1}{2}$.
\begin{proof}
This is the output of the authors' \cite{PT2} using the rigorous method described in \cite{Platt}.
\end{proof}
\end{Lem}

Finally, we need to bound $N(\sigma, T)$ which is the number of zeroes in the box $\sigma < \beta \leq 1$ and $0\leq \gamma \leq T$. We use a recent, explicit zero-density result by Kadiri, Lumley, and Ng \cite{KLN}. 

\begin{Lem}\label{sleepy}
Let $\sigma$ be a fixed number in the interval $[0.75,1)$. Then there exist positive constants $C_1(\sigma)$ and $C_2(\sigma)$ such that
\begin{equation}\label{peru}
N(\sigma, T) \leq C_1(\sigma) T^{8(1-\sigma)/3} \log^{5 - 2\sigma} T + C_2(\sigma) \log^{2}T.
\end{equation}
\end{Lem}
Table 1 of \cite{KLN} gives values for $C_1$ and $C_2$ for various\footnote{We also note that the first column in Table 1 in \cite{KLN} should have $\sigma$ in place of $\sigma_{0}$, and that the condition `$\sigma \geq \sigma_{0}$' should be deleted from the caption. We thank Habiba Kadiri for pointing this out to us. } $\sigma$, based in part on a result of the first author \cite{Platt} showing that the Riemann hypothesis holds to height $H=3.06\cdot 10^{10}$.  A reworking of \cite{KLN} using Lemma \ref{pen} would improve the constants $C_1$ and $C_2$ and hence our results. The $\log^{2}T$ term in (\ref{peru}), while  a nuisance to carry through the calculations, is utterly negligible in the final bounds. Finally, we note that Kadiri, Lumley, and Ng produce slightly superior versions of (\ref{peru}), but we have opted for the simpler version presented above.

\section{Pintz's method}\label{whim}
We follow the argument given by Pintz \cite[pp.\ 214-215]{Pintz}. We aim to derive a bound valid for all $x>x_0$. To control the $\vartheta$ term in (\ref{honduras}) we set $T=\exp(2\sqrt{\log x/R})$. We then break the sum in  (\ref{honduras}) into two pieces: those zeroes $\rho$ with $\Re \rho < \sigma$ for some $\sigma\in[1/2,1)$ to be fixed later depending on $x_0$, and the rest. For the former we have

\begin{equation}\label{bashful}
\left|\sum_{\substack{|\gamma|\leq T\\ \Re\rho<\sigma}} \frac{x^{\rho-1}}{\rho}\right| < x^{\sigma-1} \sum_{|\gamma|\leq T} \frac{1}{|\gamma|} \leq x^{\sigma-1}\left(\frac{1}{2\pi} \log \left(\frac{T}{2\pi}\right)^{2}+1.8642\right),
\end{equation}
by Lemma \ref{happy}. Note that we could improve (\ref{bashful}) by taking advantage of Lemma \ref{pen}. Such an alteration to (\ref{bashful}) makes only a negligible improvement to the bound in Theorem \ref{caprice}.
%\footnote{We could set out to replace the leading $x^{\sigma-1}$ with $\exp(-2\sqrt{\log x/R})$ for all $x>x_0$ by taking $\sigma$ sufficiently small. Specifically, we would need
%\begin{equation*}\label{dobson}
%\sigma\leq 1-\frac{2}{\sqrt{R\log x_0}}.
%\end{equation*}
%However, we will wash this term into the leading constant.}

We turn now to zeroes $\rho$ with  $\beta \geq \sigma$. We shall estimate their contribution by using the zero-free region $\sigma \geq 1 - \eta(t) = 1- (R\log t)^{-1}$ and the zero-density estimate, but by Lemma \ref{pen}, we are free to start counting from any height $H\leq 3\cdot 10^{12}$. We can therefore write
\begin{equation*}\label{latchford}
\begin{split}
\left|\sum_{\substack{|\gamma|\leq T\\ \Re\rho\geq\sigma}} \frac{x^{\rho-1}}{\rho}\right|&\leq2 \int\limits_H^T \frac{x^{-\frac{1}{R\log t}}}{t}\dif N(\sigma,t)\\
&\leq 2\left[\frac{x^{-\frac{1}{R\log T}}}{T}N(\sigma,T)-\int\limits_H^T\frac{\dif}{\dif t}\left(\frac{x^{-\frac{1}{R\log t}}}{t}\right)N(\sigma,t)\dif t\right].
\end{split}
\end{equation*}
Now the integrand only becomes negative once $t>t_0=\exp(\sqrt{\log x/R})$ (which may or may not exceed $H$) so we can write
\begin{equation*}\label{latchford1}
\begin{split}
\left|\sum_{\substack{|\gamma|\leq T\\ \Re\rho\geq\sigma}} \frac{x^{\rho-1}}{\rho}\right|&\leq 2\left[\frac{x^{-\frac{1}{R\log T}}}{T}N(\sigma,T)-\int\limits_{t_0}^T\frac{\dif}{\dif t}\left(\frac{x^{-\frac{1}{R\log t}}}{t}\right)N(\sigma,t)\dif t\right]\\
&\leq 2N(\sigma,T)\left[\frac{x^{-\frac{1}{R\log T}}}{T}+\exp\left(-2\sqrt{\frac{\log x}{R}}\right)\right]\\
&\leq \frac{2.0025 N(\sigma,T)}{T},
\end{split}
\end{equation*}
where we have used $\log x \geq 1000$ and $T=\exp(2\sqrt{\log x /R})$. We invoke Lemma \ref{sleepy} to give

\begin{equation*}\label{pejic}
\begin{split}
\left|\sum_{\substack{|\gamma|\leq T\\ \Re\rho\geq\sigma}} \frac{x^{\rho-1}}{\rho}\right|&\leq2.0025 \frac{1}{T}\left[C_1(\sigma) T^{\frac{8(1-\sigma)}{3}}\log^{5-2\sigma}T+C_2(\sigma)\log^2 T\right]\\
&\leq2.0025 \Bigg[C_1(\sigma)\exp\left(\left(\frac{10-16\sigma}{3}\right)\sqrt{\frac{\log x}{R}}\right)\left(2\sqrt{\frac{\log x}{R}}\right)^{5-2\sigma}\\
&\;\;\;\; {}+4 C_2(\sigma)\exp\left(-2\sqrt{\frac{\log x}{R}}\right)\frac{\log x}{R}\Bigg].
\end{split}
\end{equation*}

We can now write
\begin{equation*}
k(\sigma,x_0)=\left[\exp\left(\left(\frac{10-16\sigma}{3}\right)\sqrt{\frac{\log x_0}{R}}\right)\left(\sqrt{\frac{\log x_0}{R}}\right)^{5-2\sigma}\right]^{-1},
\end{equation*}
followed by
\begin{equation*}\label{manuka}
\begin{split}
C_3(\sigma,x_0)&=2\exp\left(-2\sqrt{\frac{\log x_0}{R}}\right)\log^2 x_0\;k(\sigma,x_0)\\
C_4(\sigma,x_0)&=x_0^{\sigma-1}\left(\frac{2}{\pi}\frac{\log x_0}{R}+1.8642\right)k(\sigma,x_0)\\
C_5(\sigma,x_0)&=8.01\cdot C_2(\sigma)\exp\left(-2\sqrt{\frac{\log x_0}{R}}\right)\frac{\log x_0}{R} k(\sigma,x_0)\\
A(\sigma,x_0)&=2.0025\cdot 2^{5-2\sigma}\cdot C_1(\sigma)+C_3(\sigma,x_0)+C_4(\sigma,x_0)+C_5(\sigma,x_0),
\end{split}
\end{equation*}
so that we finally reach
\begin{equation*}\label{king}
\left|\frac{\psi(x)-x}{x}\right|\leq A(\sigma,x_0)\left(\frac{\log x}{R}\right)^{\frac{5-2\sigma}{2}}\exp\left(\frac{10-16\sigma}{3}\sqrt{\frac{\log x}{R}}\right),
\end{equation*}
for all $\sigma\in[0.75,1)$ and $\log x \geq 1\,000$ large enough so that $A(\sigma,x)$ is decreasing for $x>x_0$.

It is possible to optimise the above approach further, by keeping small negative terms in play. We have not done this since, as noted after Lemma \ref{sleepy}, such an improvement ought to be coupled with a reworking of the constants $C_{1}(\sigma)$ and $C_{2}(\sigma)$ taking into account Lemma \ref{pen}. These constants may be further improved using results in \cite{CT}, and then spliced with other zero-density results in \cite{Kadiri} and \cite{Simonic}. Finally, the constant `2' in (\ref{honduras}), while not significant for Dudek's application in \cite{Dudek} plays a non-negligible role here and could well be improved. These options are all avenues for future research: we merely indicate here the power of using an explicit zero-density estimate such as (\ref{peru}).

Writing 
\begin{equation}\label{senegal}
\epsilon_0(\sigma,x_0)=A(\sigma,x_0)\left(\frac{\log x_0}{R}\right)^{\frac{5-2\sigma}{2}}\exp\left(\frac{10-16\sigma}{3}\sqrt{\frac{\log x_0}{R}}\right),
\end{equation}
we look for the $\sigma$ that results in the minimum for a given $x_0$, subject to the restriction that we need $A(\sigma,x)$ to be decreasing for all $x>x_0$. The bounds in (\ref{senegal}), having used (\ref{honduras}) are only valid when $x$ is half an odd integer. However, as an example, the difference between $\epsilon_0(0.98,\exp(1000))$ and $\epsilon_0(0.98,\lfloor\exp(1000)\rfloor -0.5)$ is less than $10^{-441}$. This approach leads to the values of $\sigma$ in Table \ref{tab3} and thus proves Theorem \ref{caprice}.

We note that it seems likely that, for larger values of $x_0$, taking $\sigma$ larger than $0.99$ would give an improvement, whence any extension to Table 1 of \cite{KLN} would be useful in this regard.

\section{Proof of Corollary \ref{white}}\label{hunger}
One can also consider the prime number theorem in the form $\pi(x) \sim \textrm{li}(x)$, where $ \textrm{li}(x) = \int_{2}^{x} (\log t)^{-1}dt$. With a little more effort we could provide a version of Corollary \ref{white} with tables of parameters, in the style of Theorem \ref{caprice}.

We begin by using partial summation and integration by parts to yield
\begin{equation}\label{green}
\pi(x) - \textrm{li}(x) = \frac{\theta(x) - x}{\log x} + \frac{2}{\log 2} + \int_{2}^{x} \frac{\theta(t) -t}{t \log^{2}t}\, dt.
\end{equation}

Let us write
\begin{equation*}\label{hall}
\int_{2}^{x} \frac{\theta(t) -t}{t \log^{2}t}\, dt = \int_{2}^{563} + \int_{563}^{e^{1000}} + \int_{e^{1000}}^{e^{2000}} + \int_{e^{2000}}^{x} = I_{1} + I_{2} + I_{3} + I_{4}.
\end{equation*}
We estimate the integral $I_{1}$ numerically. 
%This is approximately $4.63$. 
The combination of this and the $2/\log 2 $ term contributes at most 7.6 to (\ref{green}).

For $I_{2}$ we use a result of Rosser and Schoenfeld \cite{Ross}, namely that for $x\geq 563$ we have $|\theta(x) -x|\leq x/(2\log x)$. Although stronger bounds are known, we can afford to be very cavalier: the values of $I_{1}$ and $I_{2}$ are insignificant for large $x$.

For $I_{3}$ we use Theorem \ref{caprice} and (\ref{shoe}). These show that that $|\theta(x) - x|\leq 1.2\cdot 10^{-5} x$ for $x\geq e^{1000}$.  

For $I_{4}$ we follow the approach from \cite{DusartT}.  For some $\alpha$ to be determined later, define
\begin{equation*}\label{billiard}
h(t) = \frac{t \exp \{ - C\sqrt{(\log t)/R}\}}{\log^{\alpha} t}.
\end{equation*}
We want to show that 
\begin{equation}\label{dagger}
\frac{\exp \{ -C \sqrt{(\log t)/R}\}}{\log^{2-B} t}\leq h'(t), \quad (t \geq e^{2000}).
\end{equation}
Since then, if (\ref{dagger}) is true, we have
\begin{equation*}\label{mustard}
|I_{4}| \leq A_{1} R^{-B} \int_{e^{2000}}^{x} \frac{\exp \{ - C\sqrt{(\log t)/R}\}}{\log^{2-B} t}\, dt \leq A_{1} R^{-B} \int_{e^{2000}}^{x} h'(t) < A_{1} R^{-B}h(x).
\end{equation*}
Now, to show that (\ref{dagger}) is true it is sufficient to show that
\begin{equation}\label{kitchen}
\log t - (\log t)^{B-1+\alpha} - \frac{C}{2} \sqrt{ (\log t)/R} - \alpha >0, \quad (t \geq e^{2000}).
\end{equation}

We therefore end up with
\begin{equation*}\label{peacock}
|\pi(x) - \textrm{li}(x)| \leq \frac{x A_{1} (\log x)^{B-1}}{R^{B}} \exp(-C\sqrt{(\log x) / R})\left\{ 1+ \Delta\right\}, \quad (x\geq x_{0}),
\end{equation*}
where
\begin{equation}\label{study}
\begin{aligned}
\Delta =& (\log x_{0})^{1-B - \alpha}\\
&+ \frac{R^{B}\exp\left\{C\sqrt{\frac{\log x_{0}}{R}}\right\} (\log x_{0})^{1-B}}{A_{1}x_{0}} \left\{ 7.6 + \int_{563}^{e^{1000}} \frac{dt}{2 \log^{3} t} + 1.2\cdot 10^{-5} \int_{e^{1000}}^{e^{2000}} \frac{dt}{\log^{2}t}\right\}.
\end{aligned}
\end{equation}
We require that $1-B< \alpha < 2 - B$ so that the first addend in (\ref{study}) is decreasing, and so that (\ref{kitchen}) may be  satisfied.
We now take $x_{0} = e^{2000}$ and $\alpha = 0.47$. We verify that (\ref{kitchen}) is true, and we find by (\ref{study}) that $\Delta \leq 6.76$. This proves the corollary.

We briefly illustrate how one may improve upon the bound (\ref{library}), and to complete such a proof for all values of $x\geq 10$, say. We can partition the interval $[2, x]$ using the points $2, 563$ and the entries in the first column of Table \ref{tab3}. Between consecutive values in Table \ref{tab3} we use the $\epsilon_{0}$ coming from the smaller value. Indeed, we can take the smallest available $\epsilon_{0}$'s from our Table \ref{tab3}, Faber and Kadiri \cite{Faber_fix}, and the results of B\"{u}the \cite{Booth2,Booth}.  We therefore have
\begin{equation*}\label{egg1}
\begin{split}
\pi(x) - \textrm{li}(x):= E(x) &= \frac{x \epsilon_{0}(x_{0})}{\log x}+ \frac{2}{\log 2} + \int_{2}^{563} \frac{|\theta(t) -t|}{t\log^{2} t}\, dt + \frac{1}{2}\int_{563}^{e^{1000}}\frac{dt}{\log^{3} t} \\
&+ \epsilon_{0}(e^{1000}) \int_{e^{1000}}^{e^{1500}} \frac{dt}{\log^{2} t} + \cdots + \epsilon_{0}(e^{9000})\int_{e^{9000}}^{e^{9500}} \frac{dt}{\log^{2} t} + \epsilon_{0}(e^{9500})\int_{e^{9500}}^{x} \frac{dt}{\log^{2} t}.
\end{split}
\end{equation*}
Note that $E(x)$ and the right-side of (\ref{library}) are increasing in $x$. Therefore, if we wish to verify (\ref{library}) in the range $x_{0} \leq x \leq x_{1}$ we need only show that
\begin{equation*}\label{sock}
E(x_{1}) \leq 1.001\frac{x_{0}A_{1} (\log x_{0})^{B-1}}{R^{B}} \exp(-C\sqrt{(\log x_{0}) / R}).
\end{equation*}
Choosing pairs of points $(x_{0}, x_{1})$ sufficiently close together --- e.g., $x_{0} = e^{k}$ and $x_{1} = e^{k + 0.1}$ --- and populating more entries in the relevant tables of $\epsilon_{0}$'s enables us to extend Corollary \ref{white} to `moderate' values of $x$.

To extend it to `small' values of $x$, we may invoke a result by B\"{u}the \cite{Booth2}, namely, that
\begin{equation*}\label{lounge}
|\pi(x) - \textrm{li}(x)| \leq \frac{\sqrt{x}}{\log x} \left\{ 1.95 + \frac{3.9}{\log x} + \frac{19.5}{\log^{2} x}\right\}, \quad (x\leq 10^{19}).
\end{equation*}
We expect that following such a procedure, after taking advantage of the suggestions made at the end of \S \ref{whim} would give a `state of the art' version of Corollary \ref{white}.

\section{Application to Ramanujan's inequality}\label{bartolo}

Using $\pi(x) \sim \textrm{li}(x) = \int_{2}^{x} (\log t)^{-1}dt$, and integrating by parts, Ramanujan noted that
\begin{equation}\label{pitch}
\pi^{2} (x) < \frac{e x}{\log x} \pi \left(\frac{x}{e}\right),
\end{equation}
for sufficiently large $x$ --- see \cite[Ch.\ 24]{Berndt}. It is an interesting, and difficult, problem to determine the last $x$ for which (\ref{pitch}) fails. 

Dudek and Platt \cite{DPlatt} showed\footnote{We are grateful to Christian Axler who identified an error in the proof given in \cite{DPlatt}. Fortunately it was easy to fix.} that (\ref{pitch}) is true for all $x\geq \exp(9658)$. This has recently been improved by Axler \cite{Axler} to $x\geq \exp(9032)$. Dudek and Platt gave good evidence that $x = 38\,358\,837\,682 $ is the largest $x$ for which (\ref{pitch}) fails, and indeed, they showed this to be so under the Riemann hypothesis. The main obstacle in moving to an unconditional version is the size of the error term in the prime number theorem. The problem is therefore well within our wheelhouse given our result in Theorem \ref{caprice}.

As an introduction to the method given by Dudek and Platt in \cite{DPlatt} consider the expansion
\begin{equation}\label{scg}
\pi(x) = x \sum_{j=0}^{4} \frac{j!}{\log^{j+1} x} + O\left( \frac{x}{\log^{6} x}\right).
\end{equation}
Expanding, and keeping track of sufficiently many terms then gives us
\begin{equation*}\label{mcg}
\pi^{2}(x) - \frac{ex}{\log x} \pi\left(\frac{x}{e}\right) = -\frac{x^{2}}{\log^{6} x} + O\left( \frac{x^{2}} {\log^{7}}\right),
\end{equation*}
which is clearly negative for all $x$ sufficiently large. Dudek and Platt then give an explicit version of (\ref{scg}) in their Lemma 2.1. This replaces the $O(x/\log^{6}x)$ term in (\ref{scg}) by upper and lower bounds with constants $m_{a}$ and $M_{a}$ (which we give in (\ref{gabba}) and (\ref{oval})), and a criterion for $x$ being sufficiently large. Rather than duplicate the proof, we have included the critical pieces below.

We need an $a(x)$ such that
\begin{equation}\label{trunk}
 |\theta(x)-x|\log^{5} x\leq x a(x).
\end{equation}
By partial summation, and integrating $\textrm{li}(x)$ by parts once, we have, from (\ref{green}) and (\ref{trunk}), that
$$
\pi(x)\leq \frac{x}{\log x} +a(x)\frac{x}{\log^6 x}+\int\limits_2^x\frac{\dif t}{\log^2 t}+\int\limits_2^x\frac{a(t)}{\log^7 t}\dif t
$$ 
and
$$
\pi(x)\geq \frac{x}{\log x} -a(x)\frac{x}{\log^6 x}+\int\limits_2^x\frac{\dif t}{\log^2 t}-\int\limits_2^x\frac{a(t)}{\log^7 t}\dif t.
$$ 
We start with a lemma, which immediately gives us the lower bound in Theorem \ref{sidelines}.
\begin{Lem}\label{partrh}
For $x\in(599,\exp(58)]$ we have
$$
|\theta(x)-x|\leq\frac{\sqrt{x}}{8\pi}\log^2 x.
$$
\begin{proof}
This follows directly from Theorem 2 of \cite{Booth2} coupled with Lemma \ref{pen}.
\end{proof}
 \end{Lem}

We are now free to use our improved bounds. We use

\begin{equation}\label{stumps}
\frac{a(x)}{\log^5 x}=\begin{cases}
\frac{2-\log 2}{2}&2<x\leq 599\\
\frac{\log^2 x}{8\pi\sqrt{x}}& 599<x\leq \exp(58)\\
\sqrt\frac{8}{17\pi}\left(\frac{\log x}{6.455}\right)^{\frac{1}{4}}\exp\left(-\sqrt{\frac{\log x}{6.455}}\right)& \exp(58)< x<\exp(1\,169)\\
462.0\left(\frac{\log x}{5.573412}\right)^{1.52}\exp\left(-1.89\sqrt{\frac{\log x}{5.573412}}\right)&\exp(1\,169)\leq x < \exp(2\,000)\\
411.5\left(\frac{\log x}{5.573412}\right)^{1.52}\exp\left(-1.89\sqrt{\frac{\log x}{5.573412}}\right)&\exp(2\,000)\leq x < \exp(3\,000)\\
379.7\left(\frac{\log x}{5.573412}\right)^{1.52}\exp\left(-1.89\sqrt{\frac{\log x}{5.573412}}\right)&x\geq\exp(3\,000).
\end{cases}
\end{equation}
Here, the first bound is trivial, the second is Lemma \ref{partrh}, the third is from Theorem $1$ of \cite{Trudgian} and the rest are from Table \ref{caprice} adjusted as per Corollary \ref{caprice1}. We split up our bounds on $a(x)$ in this way owing to the awkward role played by the constant in the zero-free region. For some ranges of $x$, the estimate $(\log x/ R)^{1/4} \exp(-\sqrt{\log x}/R)$ is actually worse when we use a smaller (and better) zero-free region constant: see, for example, the use of $5.573412$ (from Lemma \ref{sneezy}) for $x\geq \exp(1169)$ and the use of $6.455$ from \cite{Trudgian} for smaller values of $x$. 
We could split up (\ref{stumps}) differently, but, as we shall see, the main obstacles to improving Theorem \ref{sidelines} lie in the region $x\geq \exp(3\,000)$. 
%It is for this reason that the better zero-free region in Lemma \ref{happy}, when fed into Corollary \ref{white}, improves upon what can be achieved with the earlier results in Dusart \cite[Thm 1.12]{DusartT} and Trudgian \cite[Thm 2]{Trudgian}.

Applying these improved bounds as per \cite{DPlatt} we fix $x_a>0$ so that $a(x)$ is non-increasing for all $x\geq x_a$ and write

\begin{equation*}
C_1=\frac{\log^6 x_a}{x_a}\int\limits_2^{x_a}\frac{720+a(t)}{\log^7 t}\dif t, \quad
C_2=\frac{\log^6 x_a}{x_a}\int\limits_2^{x_a}\frac{720-a(t)}{\log^7 t}\dif t, \quad C_3=2\frac{\log^6 x_a}{x_a}\sum\limits_{k=1}^5\frac{k!}{\log^{k+1} 2},
\end{equation*}
\begin{equation}\label{gabba}
M_a(x)=120+a(x)+C_1+(720+a(x_a))\left(\frac{1}{\log x_a}+\frac{7\cdot2^8}{\log^2x_a}+\frac{7\log^6 x_a}{\sqrt{x_a}\log^8 2}\right),
\end{equation}
\begin{equation}\label{oval}
m_a(x)=120-a(x)-(C_2+C_3)-a(x_a)\left(\frac{1}{\log x_a}+\frac{7\cdot2^8}{\log^2 x_a}+\frac{7\log^6 x_a}{\sqrt{x_a}\log^8 2}\right),
\end{equation}
$$
\epsilon_{M_a}(x)=72+2M_a(x)+\frac{2M_a(x)+132}{\log x}+\frac{4M_a(x)+288}{\log^2 x}+\frac{12M_a(x)+576}{\log^3 x}+\frac{48M_a(x)}{\log^4 x}+\frac{M_a(x)^2}{\log^5 x}
$$
and
$$
\epsilon_{m_a}(x)=206+m_a(x)+\frac{364}{\log x}+\frac{381}{\log^2 x}+\frac{238}{\log^3 x}+\frac{97}{\log^4 x}+\frac{30}{\log^5 x}+\frac{8}{\log^6 x}.
$$

We now need only find an $x>x_a$ such that 
$$
\epsilon_{M_a}(x)-\epsilon_{m_a}(x)<\log x.
$$
Choosing $x_a=\exp(3\,914)$ and $x=\exp(3\,915)$ will suffice, yielding the upper bound in Theorem \ref{sidelines}. 
%\end{comment}

\subsection*{Acknowledgements}
We wish to thank Habiba Kadiri, Nathan Ng, and Allysa Lumley, and the referees for very useful feedback.

 \end{document}